\documentstyle{amsppt}
\magnification=\magstep1
\vsize=20.5 true cm
\hsize=15 true cm
\voffset=0.5 true cm
\pageno=1
\nologo

\topmatter
\title
Any 3-manifold 1-dominates at most\\
 finitely many geometric 3-manifolds
\endtitle
\author
Shicheng Wang and Qing Zhou
\endauthor
\thanks{both authors are supported by Outstanding Youth Grants of NNSFC and Grants of QSSTF.}
\endthanks
\leftheadtext{S.C. Wang and Q. Zhou}
\rightheadtext{Finiteness on 1-dominating}
\endtopmatter

{\bf \S 1 Introduction.}
\vskip 0.5 true cm
Maps between 3-manifolds has been studied by many people long times ago,
and become an active subject again after
Thurston's revolution on 3-manifold theory. We refer
to [BW], [LWZ1] for various results and references on the subject.

This paper addresses
the following natural question which
was raised around 1990, see also Kirby's Problem List, [K, 3.100].

\proclaim{Question 1}
Let $M$ be a closed orientable $3$-manifold.
Are there at most finitely
many closed, irreducible and  orientable $3$-manifolds $N$ such that
there exists a degree one map $f: M\to N$?
\endproclaim

{\bf Remarks on the conditions in Question 1.}

(i) If  Poincare Conjecture  fails, i.e.,
there is a homotopy 3-sphere $N$ which is not $S^3$, then one can get
infinitely many reducible  homotopy 3-spheres by doing connected sums on $N$.
Since there always exists
degree one map from a 3-manifold $M$ to a homotopy 3-sphere, the condition
``irreducible'' on the target $N$ is posed to
avoid this unclear case.

(ii) The condition ``closed'' is posed on $M$ and $N$
just for simplicity. Indeed we can replace ``closed'' by ``compact'',
and meanwhile replace ``degree one map'' by ``degree one proper map''.
A map $f:M \to N$ between compact manifolds is proper if
$f^{-1}(\partial N)=\partial M$.

For simplicity, we adapt the following definition from [BW].
Let $M$ and $N$ be two compact orientable 3-manifolds.
Say  $M$ {\it ($1$-)dominates} $N$
if there is a proper map $f: M\to N$
of non-zero degree (degree $1$).

A closed orientable 3-manifold is called {\it geometric} if it admits one
of the following geometries:
$H^3$ (hyperbolic), $\widetilde {PSL_2 (\Bbb R)}$, $H^2\times E^1$, Sol, Nil,
$E^3$, $S^2\times E^1$, $S^3$ (spherical).
Thurston's geometrization conjecture claims that any closed, irreducible, and
orientable 3-manifolds is
either geometric or can be decomposed by
the Jaco-Shalen-Johannson torus decomposition
so that each piece is geometric.
(For details see [Th2], [Th3] or [Sc].)
All geometric 3-manifolds are precisely the
Seifert manifolds except those carry  hyperbolic geometry or  Sol geometry,
and  all geometric manifolds have infinite
fundamental group except those carry spherical geometry.

It is natural to study Question 1
when the targets are geometric first.
There are many partial results of  Question 1:

(i) The answer is affirmative if both the domain and the target
are  Seifert manifolds with infinite fundamental group [Ro], which is
based on Waldhausen's 3-manifold topology  argument.

(ii) The answer is affirmative if the domain is non-Haken
and the target is geometric [RW], which is
based on Culler-Shalen's character variety theory of 3-manifold groups.
Since the domain is non-Haken,
then the geometry of the target must be either hyperbolic or spherical.

Note also that there are  additional conditions posed on the domains
in (i) and (ii).
Two substantial result to the Question 1 are obtained recently,
where no additional conditions are posed on the domains.

(iii)  The answer is affirmative if the targets
are hyperbolic [So], which is based on the argument of
Thurston's original approach on the deformation of acylindrical manifolds.

(iv) The answer is affirmative if the target are  spherical [LWZ2],
which is based on the old knowledges of linking pair of 3-manifolds
and of combinatorial groups.

In this paper we will prove the affirmative answer to  Question 1
when the targets are all the remaining geometric 3-manifolds.
The main result of this paper is the following.

\proclaim{Theorem 1} Any orientable closed $3$-manifold
$M$ $1$-dominates at most finitely many closed orientable 3-manifolds
which are either Seifert manifolds with
infinite fundamental group or Sol manifolds.
\endproclaim

Combining Theorem 1 with  the results of [So] and [LWZ2]
we obtain the following assertion.

\proclaim{Corollary 1}
Any closed orientable $3$-manifold $1$-dominates at most finitely many
geometric $3$-manifolds.
\endproclaim

If an irreducible 3-manifold has non-trivial JSJ-decomposition, then
each decomposition piece is either a hyperbolic 3-manifold or a Seifert
manifolds with torus boundary.
From the proof of Theorem 1, we have the following corollary,
which should be useful in the discussion of non-trivial JSJ-torus
decomposition case.

\proclaim{Corollary 2} Any compact orientable $3$-manifold
$M$ dominates at most finitely many Seifert manifolds with
 non-empty boundary or zero Euler number.
\endproclaim

In Section 2, we first explain that in the proof of the main result,
one need only to deal with
Seifert manifolds with orientable orbifold base and
the torus bundle over the circle with Anosov monodromy. Then we present
various known results about degree one map, Seifert manifolds, Thurston norm
and volume of representations, including  the brief descriptions of
Thurston norm and of volume of representations.
Those results will be used in the
proof of the main results. The proof of the main result
is given Section 3.

\vskip 0.5 true cm
{\bf \S 2.  Reductions and preliminary results.}
\vskip 0.5 true cm
Each Seifert manifold has an orbifold base which is
either orientable or non-orientable.

\proclaim{Lemma 1} If there is a closed orientable 3-manifold 1-dominates
infinitely many closed orientable Seifert manifolds with
non-orientable orbifold base,
then there is a closed orientable 3-manifold 1-dominates
infinitely many closed orientable Seifert manifolds
with orientable orbifold bases.
\endproclaim

\demo{Proof} Suppose $f_j: M\to N_j$ is degree one map for all
$j\in \Bbb N$, $p_j: N_j \to O_j$ is the projection from the closed
orientable Seifert manifolds onto the non-orientable orbifold base $O_j$,
and $N_i$ and $N_j$ are not homeomorphic if $i \ne j$.
Let $\tilde q_j: \tilde O_j\to  O_j$ be the unique orientable double cover
of $O_j$ and $q_j: \tilde N_j\to N_j$ be the double covering which
covers $\tilde q_j$.
Then $\tilde N_j$  is a closed orientable Seifert manifold
with orientable orbifold base $\tilde O_j$.
Since $f_j :M\to N_j$ is of degree one, $f_{j*}: \pi_1(M)\to \pi_1(N_j)$
is onto. This implies that the index of $f^{-1}_*(\pi_1(\tilde N_j))$
in $\pi_1(M)$ is two.
Let $\tilde M_j$ be the double cover
of $M$ corresponding to the subgroup
$f_*^{-1}(\pi_1(\tilde N_j))$.
Then $f_j: M\to N_j$ can be covered by
a degree one map $\tilde f_j: \tilde M_j \to \tilde N_j$.
Since any finitely presented group has only finitely many subgroup of
given index, $\pi_1(M)$ has only finitely many subgroup of index 2.
It follows that there are only finitely many homeomorphic types
among $\{\tilde M_j; j\in \Bbb N\}$.
By passing to a subsequence, we may assume all
$\tilde M_j=\tilde M$ and we have degree one map $f_j: \tilde M\to \tilde N_j,
j\in \Bbb N$. Since any double covering is a regular covering,
each orientable Seifert manifold
double covers at most
finitely many Seifert manifolds by [MS]. It follows the homeomorphic types
of $\{\tilde N_j\}$ are infinite.  \qed

Each Sol manifold is either a tours bundle over the circle or a union
of two twisted $I$-bundle over Klein bottle.

\proclaim{Lemma 2} If there is a closed  orientable 3-manifold 1-dominates
infinitely many Sol manifold which are unions of two twisted $I$-bundle
over Klein bottle,
then there is a closed orientable 3-manifolds 1-dominates
infinitely many Sol manifolds which are torus bundle over $S^1$.
\endproclaim

\demo{Proof} Since each union of twisted $I$-bundle over Klein bottle
is double covered by a torus bundle over the circle, the rest of the
proof is the exactly
same as that we did in the proof of Lemma 1.
\qed

Let $N$ be an orientable Seifert fibered space with orientable orbifold base $F_g$
with $n$ exceptional fibers.  Then $N$ has the standard form
$(g, b;a_1,b_1;a_2,b_2;.....;a_n,b_n)$, $a_i > b_i > 0$.
There are two invariants
associated with $N$:
the Euler characteristic of the orbifold base
$$\chi_N=2-2g-\sum_{i=1}^n(1-\frac 1{a_i}),$$
and the  Euler number of $N$
$$e(N)=-b-\sum_{i=1}^n\frac {b_i}{a_i}.$$

We now give a brief description for the volume of representation
(see [Re], [Gr], [Th3] for more details).
Let $G$ be a semisimple Lie group and $X=G/K$, where $K$ is the maximum
compact subgroup of $G$. For any
orientable closed manifold $M$ and any representation $\phi: \pi_1(M)\to G$,
there is a flat $X$-bundle over $M$,
$M\times _\phi X=\frac{\tilde M\times X}{\pi_1(M)}$,
with structure group $G$, where $\tilde M$ is the universal cover of $M$,
$\pi_1(M)$ acting on the first factor $\tilde M$ by covering transformations,
and by $\phi$ on the second factor $X$.
For simplicity, we assume that $dim X = dim M=3$ and $X$ is contractible.
Let $\omega'$ be the $G$-invariant volume form on $X$, which is a closed 3-form.
Let $q: \tilde M\times X$ be the projection to the second factor. Then
$q^*(\omega')$ is a $\pi_1(M)$-invariant
closed 3-form on $\tilde M\times X$, and which induces a 3-form
$\omega$ on $M\times _\phi X$. Let $s: M\to M\times _\phi X$ be a section.
(Since $X$ is contractible, such a section exists and all such sections
are homotopic.)
We call
$\int _M s^*(\omega)\in {\Bbb R}$ the volume
of the representation $\phi$, denoted by $Vol(\phi)$, clearly
it is independent of the choice of the section $s$. Define
$$Vol_G(M)=\text{max}\{|Vol (\phi)| ; \phi :\pi_1(M)\to G\}.$$

Note if some $\phi: \pi_1(M)\to G$ is discrete and faithful,
then $M$ support the geometry of $(G,X)$ and $Vol_G(M)=Vol(\phi)$.
We get the famous Gromov norm in the case $(G, X)=(PSL_2(\Bbb C), H^3)$,
and we are interested the case $(G,X)=
(PSL_2(\Bbb R)\ltimes \Bbb Z, \widetilde {PSL_2(\Bbb R)})$ in this paper.
For short we use $SV(M)$ to denote
$Vol_{{PSL_2 (\Bbb R)\ltimes \Bbb Z}}(M)$.

\proclaim{Lemma  3}
Let $M$ and $N$ be closed orientable $3$-manifolds.
If $f: M\to N$ is a degree one map, then

(1)  $Tor H_1(M,{\Bbb Z})=A \oplus Tor H_1(N,{\Bbb Z})$,
where $Tor H_1$ is  the torsion part of $H_1$.

If $f: M\to N$ is a  map of degree $d \ne 0$, then

(2) $SV(M)\geqslant d SV(N)$.

(3) $[\pi_1(N): f_*(\pi_1(M))]|d.$
\endproclaim

\demo{Proof} For (1) see [Br, 1.2.5 Theorem].
For (2) see [BG] or [Re].
(3) is well-known and can be obtained directly
by applying covering space argument. \qed

Let $N$ be a Seifert manifold with the standard form
$$(g, b;a_1,b_1;a_2,b_2;......;a_n,b_n), \qquad a_i > b_i > 0.$$

\proclaim{Lemma 4}

(1) $N$ supporting the geometry of either $\widetilde {PSL_2 {\Bbb R}}$, or {\rm Nil},
or $H^2\times E^1$ is characterized by either
$e (N) \ne 0$ and  $\chi_N <0$, or $e(N) \ne 0$ and $\chi_N =0$, or $e(N)=0$ and $\chi_N <0$
respectively.

(2) If $e(N) \ne 0$, then the order of the torsion part of $H_1(M,\Bbb Z)$,
$$|Tor H_1(N,\Bbb Z)|=\left|\left(\prod_{i=1}^n a_i\right)
\left(b+\sum_{i=1}^n\frac {b_i}{a_i}\right)\right|
=\left|e(N) \prod_{i=1}^n a_i\right|.$$

(3) If $N$ supports the geometry of $\widetilde {PSL_2(\Bbb R)}$,
Then
$$SV(N)=\left|\frac {\chi_N^2}{e(N) }\right|.$$

(4) The equation $\chi_N=2-2g-\sum_{i=1}^n(1-\frac 1{a_i})=0$
has only finitely many solutions $(g, a_1,...,a_n)$.

(5) If $\chi_N <0$, then $\chi_N \leqslant -\frac{1}{42}$.
\endproclaim

\demo{Proof} For (1) see [Sc].
For (2) see [LWZ1, 3.1]. For (3) see [BG].
(4) and (5) are well-known and can be
obtained by elementary algebra.
\qed

Now we give a brief description on Thurston norm.
In a closed oriented 3-manifold $N$, each element
$y\in H_2(N,\Bbb Z)$ can be represented by an embedded oriented surface $F$.
Let $\chi_-(F)=\text {max}\{0,- \chi(F)\}$ if $F$ is connected, otherwise
$\chi_-(F)=\sum \chi_-F_i$, where $F_i$ are components of $F$.
Then let
$$X(y)= min\{\chi_-(F); \text{$F$ is an embedded surface representing $y$} \}.
$$
Similarly we can define $X_s(y)$ if we replace ``embedded surfaces''
by ``singular surfaces'' in the definition of $X$ (see [Th1] for details).
$X$ and $X_s$ can be extended to the second homology $H_2$ with real coefficent and are often 
called Thurston norm and Thurston singular norm respectively.

\proclaim{Lemma 5}

(1) $X$ is a pseudonorm on $H_2(M,\Bbb R)$,
in particular $mX(y)-nX(z)\leqslant X(my+nz)\leqslant mX(y)+nX(z)$.

(2) $X=X_s$.
\endproclaim

\demo{Proof} For (1) see [Th1], and for (2) see [Ga].\qed

Recall that there are only finitely many 3-manifolds support the geometries of
$S^2\times E^1$ and $E^3$, and a torus bundle over the circle
is a Sol manifold if and only if the gluing map is Anosov.
With Lemma 1 and Lemma 2, (1) of Lemma 4 to prove
Theorem 1, we need only to prove the following

\proclaim{Proposition 1}
Any orientable closed $3$-manifold
$M$ $1$-dominates at most finitely many closed orientable
3-manifolds $N_j$,
where $N_j$ belongs to one of the following classes:

(a)  Seifert manifolds with orientable orbifold bases  with
Euler number $e= 0$ and Euler characteristic $\chi< 0$.

(b)  Seifert manifolds with orientable orbifold bases  with
Euler number $e\ne 0$ and Euler  characteristic $\chi\leqslant 0$.

(c) torus bundle over the circle with Anosov monodromy.

\endproclaim

\vskip 0.5 true cm
{\bf \S 3. Proof of the Theorems.}
\vskip 0.5 true cm

In this section we will prove Proposition 1. Suppose 
contrarily that there is an orientable closed $3$-manifold
$M$ $1$-dominates infinitely many 3-manifolds $N_j$,
where $N_j$ is subject to the conditions in Theorem 1.
By passing to a subsequence we may assume that
all $N_i$'s belong to one of the following classes:

(a)  Seifert manifolds with orientable orbifold bases  with
Euler number $e= 0$ and Euler characteristic $\chi< 0$.

(b)  Seifert manifolds with orientable orbifold bases with
Euler number $e\ne 0$ and Euler characteristic $\chi\leqslant 0$.

(c) torus bundle over the circle with Anosov monodromy.

We will show that none of those three cases can happen.

Let $f_j: M\to N_j$ be a degree one map defining 1-domination.
By (1) of Lemma 3, we have
$$|Tor H_1(M,{\Bbb Z})|\geqslant |Tor H_1(N_j,\Bbb Z)|. \tag 1$$

In the first two case, we have the Seifert manifold
$$N_j=(g_j, b_j;a_{j1},b_{j1}; ......a_{jn_j},b_{jn_j}), \qquad
a_{ji} > b_{ji} > 0.$$
Since $f_{i*}:\pi_1(M)\to \pi_1(N_j)$ is surjective by (3) of Lemma 3,
the rank of $\pi_1(N_j)$ is at most the rank of $\pi_1(M)$.
The rank $\pi_1(N_j)$ is at least $2g_j+n_j-2$ [BZ], so $g_j$'s and $n_j$'s are
bounded.
Passing to a subsequence we may assume that $g_j=g$ and $n_j=n$ and we have
$$N_j=(g, b_j;a_{j1},b_{j1}; ......a_{jn},b_{jn}), \qquad a_{ji} > b_{ji} > 0\tag 2$$
Below we use $e_j$ to denote $e(N_j)$ and $\chi_j$ to denote $\chi_{N_j}$.

\vskip 0.5 true cm
\demo{Proof of Case (a)}

Each homology class $y$ of $H_2(N_j,\Bbb Z)$ can be presented
by an incompressible surface.
Since $N_j$ is irreducible Seifert manifold,
each incompressible surface is either a vertical torus
(foliated by Seifert circles), or a horizontal surfaces (transversal to
all Seifert circles) [p. 109, J].
Since $e_j=0$, $N_j$ admits horizontal surfaces.

Let $O_j$ be the orbifold base of $N_j$ and $C_j$ be a regular fiber of $N_j$.
Suppose also that $O_j$, $C_j$ and $N_j$ are compatible oriented.

Let $F_j$ be the  horizontal surface of $N_j$,
and $p_j: F_j\to O_j$ is the branched covering.
Then we have $\chi(F_j)=|d|\times \chi_j<0$,
where $d=deg(p_j)$ equals to $F_j\cap C_j$, the algebraic
intersection number of $F_j$ and a regular Seifert fiber of $N_j$.
Note that $|F_j\cap C_j|$, the absolute value of algebraic intersection number,
is precisely the geometric intersection number.

Suppose further $F_j$ is a minimal genus horizontal surface of $N_j$,
thus $F_j$ is characterized by that $|F_j\cap C_j| > 0$ is minimal.

Let $X_j$ be the Thurston norm on $H_2(N_j,\Bbb R)$.
Let $V_j = \{ y \in H_2(N_j,\Bbb Z); X_j(y) = 0 \}$,
which is generated by vertical tori.
Then $V_j$ is a subgroup of $H_2(N_j,\Bbb Z)$.

\proclaim{Lemma 6} $H_2(N_j,\Bbb Z)=\langle [F_j]\rangle +V_j$. \endproclaim

\demo{Proof} Pick any homology class $y\in H_2(N_j,\Bbb Z)$.
If $X_j(y)=0$, then $y\in V_j$. Suppose $X_j(y)\ne 0$.
Let $F$ be a oriented horizontal surface representing $y$
with $-\chi(F)=X_j(y)$.
We may assume that the degree of $p_j: F\to O_j$ is positive
(otherwise replace $y$ by $-y$).
Then $(l+1)X_j([F_j])>X_j(y)\geqslant lX_j([F_j])$  for some positive integer $l$,
that is
$$ F_j\cap C_j> (F-lF_j)\cap C_j\geqslant 0.$$
Since $F_j\cap C_j$ is minimal among all positive intersections,
we have $(F-lF_j)\cap C_j$ is zero, and therefore the minimal genus
incompressible surface which represents $[F-lF_j]$ must be a  union of
tori. That is $[F-lF_j]\in V_j$ and $y=[lF_j]+[F-lF_j]$.\qed

\proclaim{Lemma 7} If $f: M\to N$ is a  map of degree $d\ne 0$, then
 $f_* : H_*(M,\Bbb R)\to  H_* (N,\Bbb R)$ is surjective.
\endproclaim

\demo{Proof}
Recall that Poincare duality $P: H^{n-q}(N,\Bbb R)\to H_{q}(N,\Bbb R)$
is given by $z^{n-q}\to z^{n-q}\cap [N]$,
where $[N]$ ($[M]$) is the fundamental class
class of $H_n(N,\Bbb R)$ ($H_n(M,\Bbb R)$) and we also have the formula
$$f_*(f^*z^p\cap [M])=z^p\cap f_* [M] \tag 3$$
for any map $f:M\to N$.

Let $z_q\in H_q(N,\Bbb R)$. Let
$y_q=\frac 1d f^* \circ P^{-1}(z_g) \cap [M]$.
Then  by (3) we have
$$f_*(y_q)=f_*(\frac 1d f^* \circ P_N^{-1}(z_q) \cap [M])$$
$$=\frac 1d P^{-1}(z_q)\cap f_*[M]
=\frac 1d P^{-1}(z_q)\cap d[N]=z_q. \qed$$

Let $X_M$ be the Thurston norm on $H_2(M,\Bbb R)$, and $X_{sj}$ be the
Thurston singular norm on $H_2(N_j,\Bbb R)$.
Let $z_1,...,z_m$ be a basis of $H_1(M,\Bbb Z)$
and $S_i$ be a surface representing
$z_i$ with $\chi_-(S_i)=X_M(z_i)$, for $i=1,...,n$.

Let $y_i=[f_{j}(S_i)]=l_{ji} [F_j] + v_{ij}$, where $v_{ij}\in V_j$.
By (1) of Lemma 5, we have
$$X_j(y_i)=X_j(l_{ji}[F_j]+v_{ji})\leqslant l_{ij}X_j([F_j])-X_j(v_{ji})
=l_{ij}X_j([F_{j}]) \tag 4$$
Then by (2) of Lemma 5 and the definition of
Thurston singular norm, we have
$$X_j(y_i) = X_{sj}(y_i) \leqslant \chi_-(S_i)=X_M(z_i) \tag 5$$
Combine (4) and (5), we have
$$X_M(z_i) \geqslant X_j(y_i) \geqslant l_{ji}X_j([F_j]). \tag 6$$
Let $L= \text{max}\{X_M(z_i); i=1,...,m\}$. If $X_j([F_j])>L$,
then (6) implies that $l_{ji}=0$,
and therefore $y_i=[f_{j}(S_i)]=v_{ij}$. It follows
that $f_{j*}(H_2(M,\Bbb Z)) \subset V_j$. 
It contradicts Lemma 7
that $f_{j*}: H_2(M,\Bbb R)\to H_2(N_j,\Bbb R)$ is surjective.

So $L \geqslant X_j(F_j)>0$.
By passing to a subsequence we may assume that
all $X_j(F_j)$ are the same, therefore all $F_j$'s have
the same homeomorphic types, denoted by $S$.

Cutting $N_j$ along the horizontal surface $S$, we obtained an $I$-bundle
over $S$, and therefore $N_j$ can be presented as a surface bundle over $S^1$
with fiber $S$ and monodromy $g_j:S\to S$.
Since $N_j$ is a Seifert manifold, $g_j$ must be a periodic map [VI. 26., J].
However it is well-known that there are only finitely many 
periodic maps on the given surface $S$ up to conjugacy.
Since any two conjugated gluing map provide the homeomrphic
3-manifolds,
there are only finitely many homeomorphic types among all $N_j$'s.
We reach a contradiction. We have proved that Case (a) cannot happen.\qed

\vskip 0.5 true cm
\demo{Proof of Corollary 2} First note that Lemma 7 is stated for any
degree $d\ne 0$, and is still true for proper maps between manifolds
with non-empty boundaries.
Then note that for manifolds with boundary, Thurston norm was established
 and Lemma 5 is still valid ([Th1] and [Ga]).
Finally note that if $N_j$ is a Seifert manifold with boundary, then $N_j$
always contains a horizontal embedded surface.
With those three facts. The proof of Corollary 2 is the same as the proof
of Case (a). \qed

\vskip 0.5 true cm
\demo{Proof for Case (b)}

By (2) of Lemma 4, we have
$$|Tor H_1(N_j,\Bbb Z)|=
\left|\left(\prod a_{ji}\right)\left(b_j+ \sum \frac {b_{ji}}{a_{ji}}\right)\right|
=\left|e_j\prod_{i=1}^n a_{ji}\right|\tag 7$$

If all $a_{ji}$'s are uniformly bounded, then
all $b_{ji}$'s are uniformly bounded. Since we assume that all $N_j$'s are
in different homeomorphic type, we must have that $b_j$ is unbounded.

Since $\prod a_{ji}>1$ and $|\sum \frac {b_{ji}}{a_{ji}}|\leqslant n$,
we have
$$\left|\left(\prod a_{ji}\right)\left(b_j+\sum \frac {b_{ji}}{a_{ji}}\right)\right|
\geqslant |b_j-n|\tag 8$$

By (1) we have $|Tor H_1(N_j,\Bbb Z)|$ is bounded for all $j$, and
by (7) and (8)  we have  $|Tor H_1(N_j,\Bbb Z)|$ is unbounded. We reach
a contradiction.

By (4) of Lemma 4, we have ruled out the situation that $\chi=0$.
Below we assume that $\chi_j< 0$, i.e., all $N_j$
support the $\widetilde {PSL_2 \Bbb R}$ geometry.

Now we assume that some $a_{ji}$ tends to infinite as $j$ tends to infinite
up to a subsequence. Then $|\prod a_{ji}|$ tends to infinite as $j$ tends to
infinite. To be not contradicted with (1) and (7),
We must have
$$|e_j|=\left|b_j+\sum \frac {b_{ji}}{a_{ji}}\right|\to 0
\qquad\text{as $j\to \infty$}\tag 9$$

Since $\chi_j<0$, we have $\chi_j\leqslant -\frac 1{42}$ by (5) of Lemma 4, and then
$|\chi_j|\geqslant \frac 1{42}$. Then by  (5) of Lemma 4,   we have
$$SV(N_j)=\left|\frac {\chi_j^2}{e_j}\right|\geqslant \left|\frac 1{(42)^2e_j}\right|,$$
which is tends to infinite as $j$ tends to infinite.
But by (2) of Lemma 4 we have
$$SV(M)\geqslant SV(N_j)\tag 11$$
i.e. $SV(N_j)$ is uniformly bounded for all $j$.
We reach a contradiction again.
We have proved that Case (b) cannot happen.
\qed

\vskip 0.5 true cm
\demo{Proof of Case (c)}

Now $N_j$ is a torus $T$ bundle over the circle with Anosov map $g_j$,
denoted as $(T, g_j)$ sometimes.
Fix a  basis of $H_1(T,\Bbb Z)$.
Let $SL_2(\Bbb Z)$ be the group of 2 by 2  invertible integer matrices.
Let $A_j\in SL_2(\Bbb Z)$ presents $g_j$ under
  the chosen basis of $H_1(T,\Bbb Z)$.
Then $A_j$ has two  real eigenvalues of $\lambda_j$ and $\lambda_j^{-1}$,
with $|\lambda_j|>1$.

Using HHN extension one can calculate directly that
$$Tor H_1(N_j,\Bbb Z) = \frac {H_1(T,\Bbb Z)}{\langle I-A_j\rangle},$$
where $I$ is the unit of $SL_2(\Bbb Z)$, and the first Betti number of $N_j$ is 1.

Then by linear algebra we have
$$|Tor H_1(N_j,\Bbb Z)|= |I-A_j|=|(1-\lambda_j)(1-\lambda^{-1}_j)|
=|(2-(\lambda_j+\lambda_j^{-1})|\tag 12$$
By (1), $|Tor H_1(N_j,\Bbb Z)|$ is uniformly bounded, then by (12),
the absolute value of the trace $A_i$,
$|\lambda_j+\lambda_j^{-1}|$, is uniformly bounded, say by
some constant $k>0$.

Let $(T, g)$ and  $(T, g')$ be
two torus bundles over the circle with Anosov maps $g$ and $g'$.
Let $A$ and $A'$ are matrices associated
with $g$ and $g'$ under the given basis of $H_1(T, \Bbb Z)$.
If $A=BA'B^{-1}$,
for some $B\in SL_2(\Bbb Z)$, then $g=hg'h^{-1}$,
is induced by $h:T\to T$ is a homeomorphism realizing $B$.
It follows easily that $(T, g)$ and $(T, g')$ are homeomorphic.

Now a contradiction in this case will follows by the following lemma.

\proclaim{Lemma 8} There are only finitely many
conjugacy classes in $SL_2(\Bbb Z)$ representing Anosov maps
with  the absolute values of the traces are bounded by $k>0$.
\endproclaim

\demo{Proof} Let
$A=\pmatrix a & b\\c & d\endpmatrix.$ Since $A$ represents an Anosov
map, $bc\ne 0$.
Suppose $|a|> k$.
Since $|a+d|\leqslant k$, we have $|d|< 2|a|$,
and then $|ad|< 2 a^2$.
The fact $ad-bc=1$ implies that $|bc|< 2a^2 +1$. In particular,
either $|b|$ or $|c|$ is at most $\sqrt 2 |a|$.
If $|b|\leqslant \sqrt 2 a$, let
$C=\pmatrix 1 & 0\\  \pm 1 & 1\endpmatrix$,
where we chose 1 if $ab>0$ and $-1$ if $ab<0$. Then
$CAC^{-1}=\pmatrix a  \mp b & * \\ *  & * \endpmatrix$
and we can make $|a\pm b|\leqslant |a|$.
 If $|c|\leqslant \sqrt 2 a$, let
$C=\pmatrix 1 & \pm 1\\  0 & 1\endpmatrix$,
$CAC^{-1}=\pmatrix a  \mp c & * \\ *  & * \endpmatrix$
and we can make $|a\pm c|\leqslant |a|$.
This concludes that if $|a|>k$, we always can get
$A_1=\pmatrix a_1 & b_1\\c_1 & d_1\endpmatrix$
which is conjugate to $A$ in $SL_2(\Bbb Z)$ and $|a_1|<|a|$.

Therefore to prove the Lemma, we may assume that $|a|\leqslant k$. Then
similarly  since $|a+d| \leqslant k$ we have $|d|\leqslant 2k$,
and then $|ad|\leqslant 2 k^2$.
Since $ad-bc=1$, $|bc|\leqslant 2k^2 +1$. In particular all entries are bounded
by $2k^2+1$. Clearly there are only finitely many such elements in
$SL_2(\Bbb Z)$.

We have proved that Case (c) cannot happen.
\qed

We have completed the proof of  Proposition 1, and therefore
 the proof of Theorem 1.

\vskip 0.5 true cm
{\bf REFERENCES.}

[BW] M. Boileau and S. C. Wang, {\it Non-zero degree maps and surface
bundles over $S^1$},  J. Diff. Geom. {\bf 43} (1996), 789--908.

[BZ] M. Boileau and H. Zieschang, {\it Heegaard genus of
closed orientable Seifert 3-manifolds}, Invent. Math. {\bf 76} (1984)
455--468.

[BG] R. Brooks and W. Goldman, {\it Volume in Seifert space}, Duke Math. J.
Vol. 51 529-545 (1984).

[B] W. Browder, Surgery on simple-connected manifolds.
Berlin-Heidelberg-New York: springer 1972.

[Ga] D. Gabai, {\it Foliations and the topology of 3-manifolds,} J. Diff. Geom.
18 (1983) 479-536.

[Gr] M. Gromov, {\it Volume and bounded cohomology,} Publ. Math. IHES 56 5-99, (1983)

[J] W. Jaco, {\it Topology of three manifolds,} Regional Conference Series
in Mathematics. 43, AMS Math. Soc. Providence, RI 1980.

[K] R. Kirby,   {\it Problems in low-dimensional topology.}
Geometric topology, Edited by H.Kazez, AMS/IP Vol. 2.
International Press, 1997.

[H-LWZ1] C. Hayat-Legrand,  S.C. Wang, and H. Zieschang, : {\it Minimal
Seifert manifolds}. Math. Ann. {\bf 308}, 673--700 (1997)

[H-LWZ2] C. Hayat-Legrant, S.C. Wang, H. Zieschang, : {\it Any $3$-manifold $1$-dominates only
finitely Seifert manifolds with finite $\pi_1$.} Preprint.

[MS] W. Meeks and P. Scott,  {\it Finite group action on 3-manifolds,}
Invent. Math. 86, 287-346 (1986)

[ReW] A. Reid, and S.C. Wang:
{\it Non-Haken 3\--manifolds are not large with respect to mappings of
non-zero degree}. Comm. in Ann. $\&$ Geom. {\bf 7}, 105-132 (1999)

[Re] A. Resnikov, {\it Rationality of secondary classes,}
J. Diff. Geom, 43 674-692 (1996)

[Ro]  Y, Rong, {\it Maps between Seifert fibered spaces of infinite
$\pi_1$}. Pacific J. Math. {\bf 160}, 143--154 (1993)

[Sc] G. Scott, {\it The geometries of $3$-manifolds}.
Bull. London Math. Soc. {\bf 15}, 401--487 (1983)

[So]  T. Soma, {\it Non-zero degree maps onto hyperbolic $3$-manifolds.}
To appear in J. Diff. Geom.

[Th1] W. Thurston, {\it A norm for the homology of 3-manifolds.}
AMS Memoire 339  99-130 (1986).

[Th2] W. Thurston, {\it Three dimensional manifolds, Kleinian groups and hyperbolic
geometry,}  Bull. Amer. Math. Soc. Vol 6, 357-388 (1982)

[Th3] W. Thurston, Three dimensional geometry and topology, Vol 1,
Princeton University Press, NJ, 1997.

\vskip 1truecm

Peking University, Beijing, China

and

East China Normal University, Shanghai, China

\bye